\documentclass[onefignum,onetabnum]{siamart190516}

\usepackage{amsfonts}
\usepackage{graphicx}
\usepackage{epstopdf}
\usepackage{algorithmic}
\usepackage[english]{babel}
\usepackage{amsmath}
\usepackage{amsfonts}
\usepackage{stmaryrd}
\usepackage{graphicx}
\usepackage{cite}
\usepackage{subcaption}
\usepackage{hyphenat}

\usepackage[margin=1.5in]{geometry}

\usepackage{tikz}
\usetikzlibrary{calc}
\usetikzlibrary{decorations.pathreplacing}
\ifpdf
  \DeclareGraphicsExtensions{.eps,.pdf,.png,.jpg}
\else
  \DeclareGraphicsExtensions{.eps}
\fi

\newcommand{\myGlobalTransformation}[2]{\pgftransformcm{1}{0}{0.4}{0.5}{\pgfpoint{#1cm}{#2cm}}}

\newcommand{\cA}{{\cal A}}

\newcommand{\R}{\mathbb{R}}
\newcommand{\N}{\mathbb{N}}

\newsiamremark{remark}{Remark}
\newsiamremark{hypothesis}{Hypothesis}
\crefname{hypothesis}{Hypothesis}{Hypotheses}
\newsiamthm{claim}{Claim}

\newsiamthm{Def}{Definition}
\newsiamthm{Theorem}{Theorem}
\newsiamthm{Lemma}{Lemma}
\newsiamthm{Bsp}{Beispiel}
\newsiamthm{Corollary}{Corollary}
\newsiamthm{Hypothesis}{Hypothesis}

\headers{Navier-Cauchy eqn for motion estimation in dynamic imaging}{B.N. Hahn, M.-L. Kienle Garrido, C. Klingenberg and S. Warnecke}

\title{Using the Navier-Cauchy equation for motion estimation in dynamic imaging}

\author{B.N. Hahn\thanks{Department of Mathematics, University of Stuttgart 
  (\email{bernadette.hahn@imng.uni-stuttgart.de}, \email{melina-loren.kienle-garrido@imng-uni-stuttgart.de}.}
\and M.-L. Kienle-Garrido\footnotemark[1]
\and C. Klingenberg\footnotemark[2]\thanks{Department of Mathematics, University of W{\"u}rzburg
  (\email{klingen@mathematik.uni-wuerzburg.de}, \email{sandra.warnecke@mathematik.uni-wuerzburg.de}).}
\and S. Warnecke\footnotemark[2]}

\usepackage{amsopn}

\ifpdf
\hypersetup{
  pdftitle={Using the Navier-Cauchy equation for motion estimation in dynamic imaging},
  pdfauthor={B.N. Hahn, M.-L. Kienle Garrido, C. Klingenberg and S. Warnecke}
}
\fi

\begin{document}
\maketitle




\bigskip



\begin{abstract}
Tomographic image reconstruction is well understood if the specimen being studied is 	stationary during data acquisition. However, if this specimen changes its position during the measuring process, standard reconstruction techniques can lead to severe motion artefacts in the computed images. Solving a dynamic reconstruction problem therefore requires to model and incorporate suitable information on the dynamics in the reconstruction step to compensate for the motion. 

Many dynamic processes can be described by partial differential equations which thus could serve as additional information for the purpose of motion compensation. 
In this article, we consider the Navier-Cauchy equation which characterizes small elastic deformations and serves, for instance, as a simplified model for respiratory motion. 
Our goal is to provide a proof-of-concept that by incorporating the deformation fields provided by this PDE, one can reduce the respective motion artefacts in the reconstructed image.  
To this end, we solve the Navier-Cauchy equation prior to the image reconstruction step using suitable initial and boundary data. 
Then, the thus computed deformation fields are incorporated into an analytic dynamic reconstruction method to compute an image of the unknown interior structure. 
The feasibility is illustrated with numerical examples from computerized tomography.
\end{abstract}




\section{Introduction}

Tomographic imaging modalities 
represent a well-known application of the theory of inverse problems. 
If the specimen is stationary during the data collection, the reconstruction process is well understood for most imaging systems \cite{Natterer01}. A dynamic behaviour of the object during measurement, however, results in inconsistent data, and standard reconstruction techniques derived under the stationary assumption lead to severe motion artefacts in the computed images \cite{artifacts_emission,artifacts_MRI_3,shepptuningfork}. This affects medical applications, for instance due to respiratory motion, as well as non-destructive testing while imaging driven liquid fronts for oil recovery studies \cite{fluid_flow_2} or while performing elasticity experiments during the scan to determine material parameters \cite{zugversuch}.

An intuitive approach for the case of few but consistent data would be to consider a quasi-static framework. However, this is only applicable if the object motion is sufficiently slow compared to the overall data acquisition time or if the motion is periodic. 
Solving the dynamic reconstruction problem in general requires to model and incorporate dynamical prior information within the reconstruction step. For individual imaging modalities like computerized tomography, magnetic resonance imaging or positron emission tomography, several methods of this type have been proposed in the literature, based on rebinning or gating the data \cite{sijbers17,lu_mackie,gating1}, a variational formulation \cite{burger_variational,implicit_MRI,zhang,gris}, exact analytic methods \cite{exact1,exact3,hahnaffine}, iterative procedures \cite{iterative2,iterative1} or approximate inversion formulas \cite{katsevich_accurate,katsevich_local,H17a}. Further, regularization techniques developed in the general context of dynamic linear inverse problems \cite{leitao1,H14,schmitt1,schmitt2,chung} have been successfully applied to imaging problems.

The most efficient way to compensate for the dynamics is to model and incorporate the motion prior in form of a deformation map $\Phi$ which describes the trajectory of the particles in the interior of the object over time. 
In general, such deformation fields are a priori unknown and have to be extracted from the measured data. Typically, parametrized motion models are employed, i.e. only a few unknown parameters need to be estimated, either via additional measurements \cite{exact1,iterative2,motion_registration_manke,reyes} or directly from the recorded tomographic data. In computerized tomography, for instance, they can be determined by detecting traces of nodal points in the sinogram \cite{lu_mackie,H17a}. For global rotations and translations, an estimation procedure using data consistency conditions is proposed in \cite{consistency_motion_estimation}. Iterative procedures are, for example, based on edge entropy \cite{katsevich11}, or perform estimation and reconstruction step simultaneously \cite{batenburg}.

Alternatively, the dynamics can be characterized in terms of velocity fields between consecutive image frames. The intensity variations in the image sequence are then linked to the underlying velocity field by the optical flow constraint, based on  the brightness constancy assumption. Recovering both velocity fields and image frames from the measured data simultaneously requires solving non-convex optimization problems of extremely large size  \cite{burger_opticalflow,schoenlieb}. 

In this article, we pursue another approach. Many dynamic processes can be described by partial differential equations, and thus, their (numerical) solution could provide the required deformation fields. 
More precisely, we consider in the following the Navier-Cauchy equation, representing linear elasticity. In applications in radiotherapy treatment planning, the respective conservation laws are employed to model respiratory motion \cite{Werner13}. 

To reduce the overall complexity and to provide a proof-of-concept that such motion priors can compensate for the dynamics, we suggest to decouple both tasks for the study in this article. Based on the provided results, the study of the joint parameter identification problem will then be subject to future work.

In Section 2, we recall the mathematical model of dynamic imaging and present the general motion compensation strategy from \cite{hahnkienle19} in the mass preserving case which assumes that the motion is known. We then derive our elastic motion model based on conservation laws in Section 3. The respective model in particular requires prescribed initial and boundary data. Therefore, we discuss suitable choices which are feasible regarding practical applications. The numerical calculation of the deformation fields is studied in Section 4. Finally, the potential of the motion model for the purpose of motion compensation is illustrated in Section 5 at the example of computerized tomography, combining the numerically computed deformation fields with our dynamic reconstruction strategy. We conclude with an outlook to expand the suggested approach towards determining an unperturbed image and the deformation fields simultaneously via a joint parameter identification problem. 

\section{Models and reconstruction strategies in dynamic imaging}

In this section, we introduce the mathematical framework to formulate and address the problem of dynamic image reconstruction. In particular, we will consider the two-dimensional case throughout the article. Further, since the motion estimation approach via the Navier-Cauchy equation is not restricted to a particular imaging modality, we want to present the motion compensation strategy in a framework covering many different modalities. A detailed introduction can be found for instance in \cite{H14,hahnkienle19}.

We start by deriving the model of the stationary setting. To be more intuitive, we first consider the example of computerized tomography (CT). 
In CT, X-ray beams are transmitted through the specimen of interest to a detector where the intensity loss of the X-rays is recorded. Is the intensity at the source position denoted by $I_0$ and the intensity at the detector position by $I_1$, then the CT-measurement for this particular configuration is given by $\log \left(\frac{I_0}{I_1}\right)$. In particular, the radiation source needs to rotate around the object to capture information from different angles of view. Due to this rotation, the data acquisition takes a considerable amount of time. The mathematical model for this imaging process is given by the Radon transform 
\begin{equation} \label{eq:Radon} 
(\mathcal{R}h)(t,y)=\int_{\R^2} h(x)\, \delta(y-x^T\theta(t))\,\mathrm{d}x, \quad (t,y)\in [0,2\pi] \times \R,
\end{equation} 
which integrates $h$ along the straight lines $\lbrace x \in \R^2 \, : \, x^T\theta(t)=y\rbrace$, i.e. along the path of the emitted X-rays. Every source position corresponds to one time instance, in particular, the unit vector $\theta(t)=(\cos(t),\sin(t))^T$ characterizes the source position at time instance $t$, while $y$ denotes the affected detector point, and $\delta$ stands for the delta distribution. The goal is then to recover $h$, the linear attenuation coefficient of the studied specimen, from measurements $g(t,y)=(\mathcal{R}h)(t,y)$ with $(t,y)\in [0,2\pi]\times \R$. Using the Fourier transform of $\delta$, we further obtain the equivalent representation
$$
(\mathcal{R}h) (t,y)=\int_{\mathbb{R}} \int_{\R^2} (2\pi)^{-1/2}\, e^{i\sigma (y-x^T\theta(t))} \, h(x)\, \mathrm{d}x\, \mathrm{d}\sigma. 
$$

Besides CT, many imaging modalities in the stationary setting can be modeled mathematically by a linear operator which integrates the searched-for quantity along certain manifolds, for instance along circles, respectively spheres, in SONAR or photoacoustic tomography. Thus, we consider in the following a more general framework, namely model operators of type
\begin{align}\label{eq:staticOperator}
\cA h (t,y)=\int_{\mathbb{R}}\int_{\Omega_x} h(x)\, a(t,y,x) e^{i\sigma (y-H(t,x))}\,\mathrm{d}x\mathrm{d}\sigma, \quad (t,y)\in \R_T\times \Omega_y,
\end{align}
where $\Omega_x$ and $\Omega_y$ denote open subsets of $\mathbb{R}^2$ and $\mathbb{R}$, respectively, $\R_T\subset \R$ represents an open time interval covering the time required for the measuring process, $a\in C^\infty(\mathbb{R}_T\times \Omega_y\times \Omega_x)$ is a weight function and $H:\R_T\times \R^2 \rightarrow \R$ characterizes the manifold we are integrating over. 

With this observation model, we can formulate the associated inverse problem: Determine $h$ from measured data 
\begin{align} g(t,y) = \cA h (t,y) ,\quad (t,y)\in \R_T\times \Omega_y.\label{eq:staticIP} \end{align}
The component $t$ of the data variable expresses the time-dependency of the data collection process. The searched-for quantity $h$ itself, however, is independent of time, i.e. \eqref{eq:staticIP} corresponds to a \emph{static} image reconstruction problem. We refer to equation \eqref{eq:staticIP} also as \emph{static inverse problem}.

\subsection{The mathematical model of dynamic imaging}

Now, we consider the dynamic case, i.e. the investigated object changes during collection of the data and is therefore characterized by a time-dependent function $f:\R_T\times \R^2 \rightarrow \R$. For a fixed time, we abbreviate $f_t:=f(t,\cdot)$, i.e. $f_t$ represents the state of the object at time instance $t$. Then, the inverse problem of the dynamic scenario reads
\begin{align}\label{eq:dynamicIP_general}\cA^{dyn} f (t,y)= g(t,y)\end{align}
with the dynamic operator $\cA^{dyn}f(t,y):=\cA f_t (t,y)$. In particular, only measurements $g(t,\cdot)$ for a single time instance encode information about the state $f_t$, which is typically not sufficient to fully recover $f_t$. In CT, only line integrals in one particular direction would be available for the reconstruction of $f_t$, which is well known to be insufficient. Thus, additional information about the dynamic behavior need to be incorporated in order to solve dynamic inverse problems.

\begin{figure}[htp]
\begin{center}
	\begin{tikzpicture}
	\draw (1,0) arc (0:180:20pt and 30pt) ; 
	\draw (1,0) arc (0:-180:37pt) ;
	\draw (1cm-40pt,0) arc (0:-180:2pt);
	\draw (1cm-44pt,0) arc (0:180:15pt) ;
	\draw (1,-1)  node {$f_t$};
	\filldraw (0,0) circle (2pt) node[align=left,   below] {$x$};
	\draw [->] (2,0) arc (130:50:30pt); 
	\draw (2.75,0.5) node {$\Phi_t^{-1}$};
	\draw[rotate around={30:(5,0)}] (5,0) ellipse (30pt and 20pt);
	\draw (4,-1) node{$f_0$};
	\filldraw (5,-0.3) circle (2pt) node[align=left, above]{$\Phi_t^{-1} x$};
	\end{tikzpicture}
	\caption{The mapping $\Phi^{-1}_t$ correlates the state $f_t$ at time $t$ to the reference state $f_0$ at the initial time.}
	\label{fig:motionfunctionsimage}
	\end{center}
\end{figure}
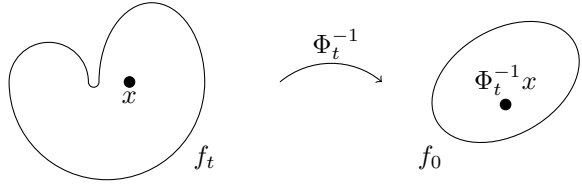

The dynamic behaviour of the object can be considered to be due to particles which change position in a fixed coordinate system of $\R^2$. This physical interpretation of object movement can then be incorporated into a mathematical model $\Phi:\R_T\times \R^2\rightarrow\R^2$, where $\Phi(0,x)=x$, i.e. we consider $f_0$ as reference state, and $\Phi(t,x)$ denotes the position at time $t$ of the particle initially at $x$. For fixed $t\in \R_T$, we write $\Phi_t x:= \Phi(t,x)$ to simplify the notation. Motivated by medical applications, where no particle is lost or added and two particles cannot move to the same position at the same time, $\Phi_t$ is assumed to be a diffeomorphism for all $t\in \R_T$. Thus, a particle $x\in \R^2$ at time $t$ is at position $\Phi_t^{-1}x$ in the reference state, see Figure \ref{fig:motionfunctionsimage}. A description of this motion model can also be found, for instance, in \cite{H14,katsevich_accurate,katsevich_local}. 

Using this motion model and the initial state function $f_0$, we find the state of the object at time instance $t$ to be
\begin{align} \label{eq:momo_masspreserving} f(t,x)=f_0(\Phi_t^{-1}x) |\det D\Phi_t^{-1}x|\end{align}
by taking into account that mass shall be preserved. \\ 

Inserting the property (\ref{eq:momo_masspreserving}) in the definition of the dynamic forward operator $\cA^{dyn}$, we obtain an operator $\cA_\Phi$ for the initial state function, namely
\begin{align}
\cA_\Phi f_0 (t,y) := \cA(\vert\det \text{D}\Phi_t^{-1}(\cdot)\vert (f_0\circ \Phi_t^{-1}))(t,y).
\label{eq:dynamicOperator}
\end{align}

\begin{remark}
In our previous work \cite{H14,H17a,hahnquinto2016}, we considered the intensity preserving model
\[f(t,x)=f_0(\Phi_t^{-1}x),\]
i.e. each particle keeps its initial intensity over time. Although this does not alter the nature of our reconstruction algorithm, we insist here on the mass preserving case to be consistent with the conservation laws employed in Section 3 for the purpose of motion estimation and clinical applications. The mass preserving model is also considered, for instance, in \cite{katsevich_accurate,katsevich_local}. 
\end{remark}

For a theoretical analysis, the motion model $\Phi$ is typically assumed to satisfy the following additional conditions, cf. \cite{siamak_rabieniaharatbar,ChungNguyen-MotionPAT,hahnquinto2016,hahnkienlequinto2020}:
\begin{itemize}
	\item The map \begin{equation}
	\label{bolker1} x \mapsto \begin{pmatrix}
	H(t,\Phi_t x)\\ D_t H(t,\Phi_tx)
	\end{pmatrix}
	\end{equation}
	is one-to-one for each $t$.
	\item It holds
	\begin{equation}\label{bolker2}
	\det \begin{pmatrix}
	D_x H(t,\Phi_tx)\\D_x D_t H(t,\Phi_t x) 
	\end{pmatrix}\neq 0
	\end{equation}
	for all $x\in\R^2$ and all $t\in \R_T$.
\end{itemize}
Basically, these properties ensure that the object's motion does not result in trivial sampling schemes for $f_0$. A detailed interpretation of these conditions can be found, for instance, in \cite{hahnquinto2016}.

If the deformation fields $\Phi_t$ are known, the dynamic inverse problem (\ref{eq:dynamicIP_general}) reduces to determining $f_0$ from the equation
\begin{align} \label{eq:dynamicIP_motion}
\cA_\Phi f_0=g.
\end{align}
In \cite{hahnkienle19,H14,katsevich_accurate}, efficient algorithms have been developed to solve this task. The underlying strategy proposed in \cite{hahnkienle19} is summarized in the following, before 
we introduce our PDE-based approach to determine the deformation fields $\Phi_t$ in Section 3 and combine both strategies to solve \eqref{eq:dynamicIP_motion} when $\Phi_t$ are unknown.

\subsection{Motion compensation algorithms}

Throughout this section, we assume the motion $\Phi$ to be known and focus on solving (\ref{eq:dynamicIP_motion}). Under suitable assumptions on the phase function $H$, the linear integral operator $\cA$ from the underlying static case belongs to the class of \emph{Fourier integral operators}. To define this type of operators, we first introduce the concepts of amplitude and phase function. 

\begin{definition}\quad\vspace{-3ex}\\
	
	\begin{itemize}
		\item Let $\Lambda\in C^\infty(\R_T\times \Omega_y\times \Omega_x\times \R\backslash \{0\})$ be a real-valued function with the following properties:
		\begin{enumerate}
			\item $\Lambda$ is positive homogeneous of degree $1$ in $\sigma$, i.e. $\Lambda(t,y,x,r\sigma)=\newline r\Lambda(t,y,x,\sigma)$ for every $r>0$,
			\item both $(\partial_{(t,y)}\Lambda,\partial_\sigma \Lambda)$ and $(\partial_x\Lambda, \partial_\sigma\Lambda)$ do not vanish for all $(t,y,x,\sigma)\in\R_T\times \Omega_y\times \Omega_x\times \R\backslash\{0\}$,
			\item it holds $\partial_{(t,y,x)}\left(\frac{ \partial\Lambda}{\partial\sigma}\right)\neq 0$ on the zero set 
			\begin{align*}
			\Sigma_\Lambda=\{(t,y,x,\sigma)\in\R_T\times \Omega_y\times \Omega_x\times\R\backslash\{0\}: \partial_\sigma \Lambda =0\}.
			\end{align*}
		\end{enumerate}
		Then, $\Lambda$ is called a \emph{non-degenerate phase function}. 
		\item Let $a\in C^\infty(\R_T\times \Omega_y\times \Omega_x\times \R)$ satisfy the following property:\\
		For every compact set $K\subset \R_T\times \Omega_y\times \Omega_x$ and for every $M\in\mathbb{N}$, there exists a $C=C(K,M)\in\R$ such that
		\begin{align*}
		\left\vert\frac{\partial^{n_1}}{\partial t^{n_1}}\frac{\partial^{n_2}}{\partial y^{n_2}}\frac{\partial^{n_3}}{\partial x_1^{n_3}}\frac{\partial^{n_4}}{\partial x_2^{n_4}}\frac{\partial^{m}}{\partial \sigma^{m}} a(t,y,x,\sigma)\right\vert\leq C(1+\vert\sigma\vert)^{k-m}
		\end{align*}
		for $n_1+n_2+n_3+n_4\leq M$, $m\leq M$, for all $(t,y,x)\in K$ and for all $\sigma\in\R$.\\
		Then $a$ is called an \emph{amplitude} (of order $k$).
		\item Let $\Lambda$ denote a non-degenerate phase function and let $a$ be an amplitude (of order $k$). Then, the operator $\mathcal{T}$ defined by
		\begin{align*}
		\mathcal{T}u(t,y)=\int u(x) a(t,y,x,\sigma) e^{i  \Lambda(t,y,x,\sigma)} \mathrm{d}x\mathrm{d}\sigma, \quad (t,y)\in \R_T\times \Omega_y
		\end{align*}
		is called a \emph{Fourier integral operator} (FIO) (of order $k-1/2$).
	\end{itemize}
\end{definition}
For more details and a more general definition see \cite{hormander2009,Treves:1980vf}.\\

In \cite{hahnkienle19,hahnkienlequinto2020}, it was shown that under suitable smoothness conditions on $\Phi$, the dynamic operator $\cA_\Phi$ inherits the FIO property from its static counterpart $\cA$.

\begin{theorem}
	Let $\Phi\in C^\infty(\R_T\times \R^2)$ and let $\Phi_t$ be a diffeomorphism for every $t\in\R_T$. If the static operator $\cA$ from \eqref{eq:staticOperator} is an FIO, the respective dynamic operator $\cA_\Phi$ from (\ref{eq:dynamicOperator}) is an FIO as well. 
\end{theorem}

Fourier integral operators have specific properties that can be used to design efficient motion compensation strategies: They  encode characteristic features of the object - the so-called \emph{singularities} - in precise and well-understood ways. 

\begin{figure}[htp]
	\begin{center}
		\includegraphics[width=4cm]{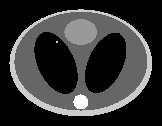}
		\hspace{2ex}
		\includegraphics[width=4cm]{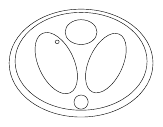}
	\end{center}
	\caption{Initial state $f_0$ of a phantom (left) and its singularities (right).}\label{fig:contours}
\end{figure}

Formally, singularities of a (generalized) function $h$ correspond to the elements of the \emph{singular support} $\text{ssupp}(h)$, which denotes the complement of the largest open set on which $h$ is smooth. In imaging applications, where the searched-for quantity is typically piecewise constant (each value characterizing a particular material), the singularities correspond to the contours of $h$, see Figure \ref{fig:contours}.

The method for motion compensation from \cite{hahnkienle19} is motivated by results on microlocal analysis, which address - among others - the question which singularities can be stably recovered from the data. The main idea is to use reconstruction operators of the form 
\begin{align}\label{eq:reconstrOperator}
	\mathcal{L}_\Phi=\mathcal{B}_\Phi \mathcal{P} 
\end{align}
on the data $g=\mathcal{A}_\Phi f_0$ with $\mathcal{P}$ a \emph{pseudodifferential operator} (typically acting on the spatial data variable $y$) and a \emph{backprojection operator} $\mathcal{B}_\Phi$ which incorporates the information on the dynamic behavior. 

\begin{definition}
	\begin{enumerate}
		\item[a)] An operator of the form
		$$
		\mathcal{P}g(t,s) = \int_\R \int_{\R} e^{i\sigma(s-y)}\, p(s,y,\sigma)\, g(t,y)\,\mathrm{d}y\,\mathrm{d}\sigma
		$$
		with $|\sigma|\leq 1$ and amplitude $p$ which is locally integrable for $s, y$ in any compact set $K$ is called \emph{pseudodifferential operator (PSIDO)} (acting on the spatial data variable $y$). 
		\item[b)] The operator 
		\begin{align*}
		\mathcal{B}_\Phi g (x) = \int_{\R_T} b(t,x)\, g(t, H(t,\Phi_t x))\,\mathrm{d}t, \quad x\in \R^2,
		\end{align*}
		where $b(t,x)$ is a positive $C^\infty$-weight function on $\R_T\times \R^2$, is called \emph{backprojection operator} associated to $\cA_\Phi$. 
	\end{enumerate}
\end{definition}

With these representations of $\mathcal{B}_\Phi$ and $\mathcal{P}$, the operator $\mathcal{L}_\Phi$ from \eqref{eq:reconstrOperator} reads
\begin{equation}\label{eq:reconstrOperatorDetail}
\mathcal{L}_\Phi g (x) =  \int_{\R_T} \int_{\R} \int_{\R} b(t,x) \, p(H(t,\Phi_tx),y,\sigma) \, g(t,y) e^{i\sigma (H(t,\Phi_tx)-y)}\,\mathrm{d}y\,\mathrm{d}\sigma \mathrm{d}t.
\end{equation}

\begin{remark}
\begin{enumerate}
	\item[a)] Pseudodifferential operators constitute a special case of an FIO. A more general definition than the one given above can be found, for instance, in \cite{KrishnanQuinto2014}.
	\item[b)] If we choose the weight $ b(t,x) = a(t,H(t,\Phi_t), \Phi_t x)$ with the amplitude $a$ of the underlying static operator $\cA$, the respective backprojection operator $\mathcal{B}_\Phi$ corresponds to the dual operator of $\cA_\Phi$.
\end{enumerate}
\end{remark}

The following result forms the basis to our motion compensation method.

\begin{theorem}
	Let $\Phi\in C^\infty(\R_T\times \R^2)$ and let $\Phi_t,\, t\in \R_T$ be diffeomorphisms that satisfy the conditions \eqref{bolker1} and \eqref{bolker2}. Further, let $\mathcal{L}_\Phi = \mathcal{B}_\Phi \mathcal{P}$ be well-defined. Then, $\mathcal{L}_\Phi$ preserves the singularities of $f_0$ which are ascertained in the measured data. 
\end{theorem}

\begin{proof}
Since $\Phi$ is smooth and $\Phi_t$ is a diffeomorphism for each $t\in\R_T$ and since $\mathcal{L}_\Phi = \mathcal{B}_\Phi \mathcal{P}$ is well-defined for the PSIDO $\mathcal{P}$, Hypothesis 1 from \cite{hahnkienlequinto2020} is fulfilled. Further conditions \eqref{bolker1} and \eqref{bolker2} are satisfied, which are the conditions (17) and (18) in \cite{hahnkienlequinto2020}. Hence, Theorem 13 from \cite{hahnkienlequinto2020} can be applied. This gives us
\begin{align*}
    \text{ssup}(\mathcal{L}_\Phi f_0) \subset \text{ssup}(f_0) \cap \mathcal{V},
\end{align*}
where $\text{ssup}(h)$ is the singular support of $h$, which is the set of singularities of $h$, for more information see \cite{hahnkienlequinto2020}. Further $\mathcal{V}$ is the set of visible singularities of $f_0$ or in other words of singularities of $f_0$ that are ascertained in the measured data.
\end{proof}

\textbf{Interpretation:} Applying a reconstruction operator $\mathcal{L}_\Phi$ of type \eqref{eq:reconstrOperator}   provides an image showing the singularities of $f_0$ correctly, which are encoded by the dynamic data. In particular, no motion artefacts arise. Thus, the described approach provides in fact a motion compensation strategy. In particular, it can be easily implemented and the computational effort is comparable to the one of static reconstruction algorithms of type \emph{filtered backprojection}. If an inversion formula of type $u = \cA^* \mathcal{P}^{stat} \cA u$ with a PSIDO $\mathcal{P}^{stat}$ is known for the static case, then choosing the PSIDO $\mathcal{P}=\mathcal{P}^{stat}$ for the motion compensation strategy provides even a good approximation to the exact density values of $f_0$ \cite{hahnkienle19}. In computerized tomography, such an inversion formula is known with $\mathcal{P}^{stat}$ being the \emph{Riesz potential} \cite{Natterer86}. \\

\begin{remark}
Although the ascertained singularities of $f_0$ are correctly reconstructed by $\mathcal{L}_\Phi$, some additional artefacts might occur if the motion is non-periodic. This has been studied in detail for computerized tomography in \cite{hahnquinto2016} and for a more general class of imaging problems in \cite{hahnkienlequinto2020}. These artefacts would be caused by singularities encoded at beginning and end of the scanning and would spread along the respective integration curve. Nevertheless, this is an intrinsic property due to the nature of the dynamic problem and therefore does not impose a major restriction to our reconstruction approach. In particular, for periodic motion as in medical applications, such as respiratory motion, the data acquisition protocol could be adjusted to the breathing cycle to avoid this issue.
\end{remark}

Since inverse problems are typically ill-posed, a regularization is required to determine $\mathcal{L}_\Phi g$ stably from the measured data $g=\mathcal{A}_\Phi f_0$. For our considered class of imaging problems, the ill-posedness is typically revealed by the growth of the symbol $p$ in terms of $\sigma$. For instance, the amplitude of the Riesz potential arising in computerized tomography corresponds to $p(s,y,\sigma) = p(\sigma) = |\sigma|$, thus, amplifying the high frequencies of the data $g$. The inversion process can be stabilized by introducing a smooth low-pass filter $e^\gamma$, i.e. by considering 
\begin{equation}\label{eq:reconstrOperatorRegu}
\mathcal{L}^\gamma_\Phi g (x) =  \int_{\R_T} \int_{\R} \int_{\R} b(t,x) \, p(H(t,\Phi_tx),y,\sigma)\, e^\gamma(\sigma) \, g(t,y) e^{i\sigma (H(t,\Phi_tx)-y)}\,\mathrm{d}y\,\mathrm{d}\sigma \mathrm{d}t
\end{equation}
with $\gamma > 0$ instead of \eqref{eq:reconstrOperatorDetail}, see \cite{hahnkienle19} for more details. 

\subsection{Reconstruction operator in dynamic CT}\label{Sect:CompensationCT}

Since we will evaluate our motion estimation strategy in Section 5 at the example of computerized tomography, we want to state the respective motion compensation algorithm for this application explicitely.

As introduced in the beginning of this section, the mathematical model operator $\mathcal{A}$ of the static case corresponds to the classical Radon transform $\mathcal{R}$, see \eqref{eq:Radon}, which is an FIO with amplitude $a(t,y,x) = (2\pi)^{-1/2}$ and phase function $\Lambda(t,y,x,\sigma) = \sigma (y-H(t,x))$, where $H(t,x) = x^T\theta(t)$ \cite{KrishnanQuinto2014}. Thus, the associated dynamic backprojection operator $\mathcal{B}_\Phi$ with weight $b(t,x)=a(t,H(t,\Phi_t),\Phi_tx) = (2\pi)^{-1/2}$ reads
$$
\mathcal{B}_\Phi g(x)= (2\pi)^{-1/2} \, \int_{\R_T} g(t, (\Phi_t x)^T\theta(t))\, \mathrm{d}t. 
$$
Choosing as PSIDO the Riesz potential with amplitude $p(s,y,\sigma) = |\sigma|$ and a low-pass filter $e^\gamma$, for instance the Gaussian, we obtain the dynamic reconstruction operator
\begin{equation*}\label{eq:reconstrOperatorDynCT}
\mathcal{L}^\gamma_\Phi g (x) = (2\pi)^{-1/2} \int_{\R_T} \int_\R \int_{\R} |\sigma|\, e^\gamma(\sigma) \, g(t,y) \, e^{i\sigma ((\Phi_tx)^T \theta(t) - y)}\,\mathrm{d}y\,\mathrm{d}\sigma \mathrm{d}t, \quad \gamma > 0,
\end{equation*}
which can be implemented in form of a \emph{filtered backprojection} type algorithm, see \cite{H17a}.

\section{Linear elastics}

In this section and the following one, we will treat the task of motion estimation. While, for a global deformation, the dynamic behavior of the boundary can be observed externally, the deformation in the interior is a priori unknown. Since many dynamic processes can be mathematically described in terms of a partial differential equation (PDE), we propose to determine the deformation fields $\Phi_t$ by finding the solution of an appropriate PDE with suitable given initial and boundary data. \\
\\
Since the deformation fields $\Phi_t, t\in \R_T$ describe the mapping from the initial/reference state to the current position, we choose the Lagrangian description for the PDE. Let $\Omega_x \subset \mathbb{R}^2$ denote the initial domain, i.e. $\Omega_x$ corresponds to the support of the initial state $f_0$, and consequently, we choose $\Omega_x$ to be the reference configuration. 

We require that $\Phi_t, t\in \R_T$ preserves its orientation meaning that $\det D\Phi(t,x)>0$ for all $(t,x)\in \R_T\times \Omega_x$. Especially in medical applications, this assumption is sensible since it also states that the local ratio of the current and the initial volume never vanishes. \cite{Antman04}\\
\\
The following definition links the current and the initial position. 
\begin{definition}
The difference between the current and the initial position is called displacement $u(t,x) = \Phi(t,x) - x$ for all $(t,x) \in \mathbb{R}_T \times \Omega_x$.
\end{definition}

Our investigations are driven by medical applications. Having the cross section of a thorax in mind, we consider two spatial dimensions, which is reasonable under a plane strain assumption. The properties of respiratory motion shall then be reflected by adequate equations. Due to its periodic behavior, it is clear that occurring stresses do not cause any yielding.  So we assume a linear relationship between stresses and strain which results in linear elasticity.  In future work, we plan to consider more general stress-strain laws.

We consider this paper a proof-of-concept. Thus we insert Hooke's law in the general equation of conservation of momentum and arrive at the Navier-Cauchy equations in two spatial dimensions for $(t,x) \in \mathbb{R}_T \times \Omega_x$, see for reference \cite{Temam_05}:
\begin{align} \label{eq:NavierCauchy}
\hat{\rho} \, \frac{\partial^2 u_k }{\partial t^2} = \hat{v}_k +  \mu \, \left( \frac{\partial^2 u_k}{\partial x_1^2} + \frac{\partial^2 u_k}{\partial x_2^2} \right) + (\lambda + \mu) \, \frac{\partial}{\partial x_k } \left( \frac{\partial u_1}{\partial x_1} + \frac{\partial u_2}{\partial x_2} \right) \quad \text{for} \quad k = 1,2.
\end{align}

These are two linear PDEs for the two unknown components $u_1, u_2$ of the displacement $u$ with the following parameters:
\begin{itemize}
\item The density $\hat{\rho}  = \rho(t,x) \operatorname{det} D\Phi(t,x) $ equals the initial density distribution \newline $\hat{\rho}= \hat{\rho}(x) = \rho(0,x)$ due to the conservation of mass.  
\item The external volume forces are denoted by $\hat{v} = v(t,x) \operatorname{det} D\Phi(t,x)$, where \newline $v: \mathbb{R}_T \times \Omega_x \to \mathbb{R}^2$ describes the volume force density.
\item The Lam\'{e}-coefficients $\lambda$ and $\mu$ specify the behavior of the material.
\end{itemize}

For a fully determined problem, we need the displacements at time $t = 0$ and their time derivatives as initial data
\begin{align*}
u(0,x) = \vartheta^0(x)  \quad \text{and} \quad
\frac{\partial }{\partial t} u(0,x) = \vartheta^1(x),
\end{align*} 
with some given $\vartheta^0, \vartheta^1: \Omega_x \to \mathbb{R}^2$.\\
\\
Also the behavior of the boundary needs to be known, more precisely a function \newline $\psi: \mathbb{R}_T \times \Omega_x \to \mathbb{R}^2$ prescribing the evolution of the displacements on the boundary of the domain $\Gamma = \partial \Omega_x$: 
\begin{align*}
u(t,x) = \psi(t,x) \quad \text{for} \quad (t,x) \in \mathbb{R}_T \times \Gamma.
\end{align*}

Solving the PDE we have introduced with given initial and boundary conditions corresponds to determining the displacement $u$, respectively the deformation $\Phi$ in the interior of the object from observations of the dynamic behavior of the object's boundary.  This way we model the movement in the object's interior, which provides exactly the information about the motion needed for our motion compensation algorithm.  \\

Under some regularity assumptions, existence and uniqueness of the solutions of the Navier-Cauchy equation (\ref{eq:NavierCauchy}) can be proven. If the initial data is $C^\infty$, solutions for the initial value problem stay $C^\infty$, cf. \cite{HughesKatoMarsden}.
Also for the initial-boundary value problem, there are existence and uniqueness results, cf. \cite{ChenvonWahl}. For appropriate boundary data $\psi$, regularity of the solutions does not get lost, and it can be shown that the solutions are diffeomorphisms, cf. \cite{Ciarlet88}. In our numerical experiments in Section 5, the initial and boundary data is chosen so that the application of the motion compensation algorithm goes through. \\

In the following, we quickly discuss suitable initial and boundary data regarding our application in dynamic imaging. As mentioned before, a global motion can be observed externally, thus, we make the reasonable assumption that the boundary data $\psi(t,x), \ (t,x)\in \R_T\times \Gamma$ are given. However, in practice, only discrete boundary data $\psi(t_n, x_{i,j})$, $n=1,\dots,N$, $i=1,\dots,I$, $j=1,\dots,J$, $N,I,J\in \N$ will be available which might be even sparse with respect to the spatial component (i.e. $I,J$ might be small) or corrupted by noise. This will be addressed in our numerical study in Section 5.

Since we are overall interested in a reconstruction of the initial state of the object and since we start with an undeformed configuration, the initial displacement data $\vartheta^0$ and $\vartheta^1$ will be set to zero. 

\begin{remark} \label{remark:initial-density}
According to \eqref{eq:NavierCauchy}, the Navier-Cauchy contains the initial density distribution $\hat{\rho}$ as parameter which is strongly linked to the quantity $f_0$ we would like to determine by our imaging modality (in particular, they share the same singularities). If we knew this parameter $\hat{\rho}$, we would already have full knowledge about the interior structure of the studied specimen. Thus, we cannot assume to know $\hat{\rho}$. Formally, we could formulate a joint motion estimation and image reconstruction approach, where we identify the parameter $\hat{\rho}$ of the PDE using the measurements from our imaging modality.  
However, to simplify the task for our proof-of-concept study, we propose another approach. In order to decouple the tasks of motion estimation via the Navier-Cauchy equation and dynamic image reconstruction, we use for the solution of the PDE a simplified prior instead of the exact density distribution $\hat{\rho}$. This is discussed in more detail in Section 5. 
\end{remark}

\begin{remark}
In this paper, we use the Navier-Cauchy equation to approximate the moving body. 
Using data at the boundary we deduce the motion of the whole body. For this we need to fix parameters of the Navier Cauchy PDE plus fix initial data. This is done by an informed guess (cf. the discussions in Section 4 and 5). Even with these approximations we find promising results in our numerical experiments, see Section 5. In future work, more elaborated PDEs are going to be considered in order to capture further details of the body's motion.
\end{remark}

\section{Numerical solution of the Navier-Cauchy equation}
We divide the given time period $t \in \R_T$ into equidistant intervals and call the time steps $t_n = n \cdot \Delta t$. We choose a Cartesian grid (not necessarily uniform) so that the discrete boundary lies on the continuous boundary, see Figure \ref{fig:interpolation}. Using central finite differences of second order for the discretization of the Navier-Cauchy equation (\ref{eq:NavierCauchy}), we obtain an explicit numerical scheme. We have chosen finite difference for our proof-of-concept study. For future studies, we plan to use a more elaborated numerical method.

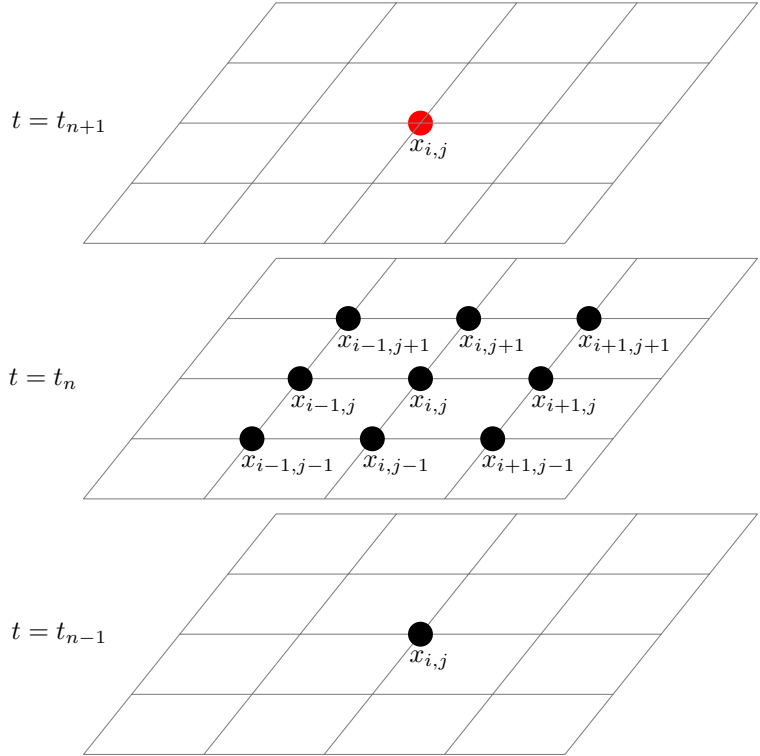
\begin{figure}[htp]
\begin{center}
\begin{tikzpicture}
\begin{scope}[scale=0.8]
        \myGlobalTransformation{0}{8.5};
        \node at (4,4) [circle,fill=red] {};
        \node at (4.25,3.75) [align=right,below]{$x_{i,j}$};
        \node at (-2,4) {$t=t_{n+1}$};
\end{scope}
\begin{scope}[scale=0.8]
        \myGlobalTransformation{0}{8.5};
        \draw [black!50,step=2cm] grid (8,8);
 \end{scope}
 
\begin{scope}[scale=0.8]
        \myGlobalTransformation{0}{0};
        \draw [black!50,step=2cm] grid (8,8);
 \end{scope}
 \begin{scope}[scale=0.8]
        \myGlobalTransformation{0}{4.25};
        \draw [black!50,step=2cm] grid (8,8);
 \end{scope}
 \begin{scope}[scale =0.8]
        \myGlobalTransformation{0}{4.25};
 				\foreach \x in {2,4,6} {
 						\foreach \y in {2,4,6} {
 								\node at (\x,\y) [circle,fill=black] {};
           			}
				}
		\node at (2.7,1.75) [align=right,below]{$x_{i-1,j-1}$};
		\node at (4.5,1.75) [align=right,below]{$x_{i,j-1}$};
		\node at (6.7,1.75) [align=right,below]{$x_{i+1,j-1}$};
		\node at (2.5,3.75) [align=right,below]{$x_{i-1,j}$};
		\node at (4.25,3.75) [align=right,below]{$x_{i,j}$};
		\node at (6.5,3.75) [align=right,below]{$x_{i+1,j}$};
		\node at (2.7,5.75) [align=right,below]{$x_{i-1,j+1}$};
		\node at (4.5,5.75) [align=right,below]{$x_{i,j+1}$};
		\node at (6.7,5.75) [align=right,below]{$x_{i+1,j+1}$};
		\node at (-2.27,4) {$t=t_n$};		
\end{scope}
\begin{scope}[scale=0.8]
        \myGlobalTransformation{0}{0};
        \node at (4,4) [circle,fill=black] {};
        \node at (4.25,3.75) [align=right,below]{$x_{i,j}$};
        \node at (-2,4) {$t=t_{n-1}$};
\end{scope}
\end{tikzpicture}
 \caption{We illustrate the stencil for our numerical scheme. For the update of the values at node $x_{i,j}$ from $t_n\to t_{n+1}$, we have to provide information about the values at the other marked nodes.}
	\label{fig:stencil}
	\end{center}
\end{figure}

We denote $x_{i,j}= ((x_1)_i,(x_2)_j)=(x_i,y_j)$, $(u_k)^n_{i,j} = u_k(t_n,x_{i,j}) $ for $k=1,2$, $\rho^0_{i,j} = \hat{\rho}(x_{i,j}) $, $\hat{v}^n_{i,j} = \hat{v}(t_n,x_{i,j})$, $\Delta x_i = x_{i+1}-x_i$ and $\Delta y_j = y_{j+1}-y_j$. Then the scheme reads exemplary for the first component $k=1$
\begin{equation*} 
\begin{split}
&(u_1)_{i,j}^{n+1} = {\scriptstyle \frac{\Delta t^2}{\rho^0_{i,j}}} \hat{v}^n_{i,j} - (u_1)^{n-1}_{i,j} + 2\left[ 1 - {\scriptstyle\frac{2\Delta t^2}{\rho^0_{i,j}} \left( \frac{\mu}{\Delta y_j^2 + \Delta y_{j-1}^2} + \frac{\lambda + 2\mu}{\Delta x_i^2 + \Delta x_{i-1}^2}\right)} \right] (u_1)^n_{i,j}\\
 &+{\scriptstyle\frac{ \Delta t^2}{\rho^0_{i,j}} \frac{2(\lambda + 2\mu)}{\Delta x_i^2 + \Delta x_{i-1}^2}} \left[
 \left( 1 - {\scriptstyle\frac{\Delta x_i - \Delta x_{i-1}}{\Delta x_i + \Delta x_{i-1}}} \right) (u_1)_{i+1,j}^{\text{n}}  + \left( 1+ {\scriptstyle\frac{\Delta x_i - \Delta x_{i-1}}{\Delta x_i + \Delta x_{i-1}}} \right) (u_1)_{i-1,j}^{\text{n}} \right] \\
&+{\scriptstyle\frac{ \Delta t^2}{\rho^0_{i,j}} \frac{2\mu}{\Delta y_j^2 + \Delta y_{j-1}^2}} \left[
 \left( 1 - {\scriptstyle\frac{\Delta y_j - \Delta y_{j-1}}{\Delta y_j + \Delta y_{j-1}}} \right) (u_1)_{i,j+1}^{\text{n}}  + \left( 1+ {\scriptstyle\frac{\Delta y_j - \Delta y_{j-1}}{\Delta y_j + \Delta y_{j-1}}} \right) (u_1)_{i,j-1}^{\text{n}} \right] \\
&+ {\scriptstyle \frac{\Delta t^2}{\rho^0_{i,j}}\frac{\lambda + \mu}{(\Delta x_i + \Delta x_{i-1})(\Delta y_j + \Delta y_{j-1})}} \left( (u_2)^n_{i+1,j+1} - (u_2)^n_{i-1,j+1} - (u_2)^n_{i+1,j-1} + (u_2)^n_{i-1,j-1} \right).
\end{split}
\end{equation*}
The corresponding stencil is illustrated in Figure \ref{fig:stencil}.

For the first time step, the (discrete) initial condition needs to be inserted
\begin{align*}
(u_k)_{i,j}^{-1} = (u_k)_{i,j}^1 - 2 \Delta t \;\vartheta^1(x_{i,j}) \quad \text{for} \quad k=1,2.
\end{align*}

The stencil for the spatial discretization has nine nodes. Since we are inspired by medical applications and a thorax is a possible specimen to be studied, we might deal with curved domains. For curved domains at the boundary, for the update scheme there is a node, which is not available to the stencil, see Figure \ref{fig:interpolation}. Hence, we need to use an interpolation method. 

For reasons of stability, we want to maintain the stencil. We call the missing node a ghost node that needs to have a value assigned to it, and we denote $h$ the quantities given at every node. The indices of the nodes are given in Figure \ref{fig:interpolation}. A second-order approach is the following one for the components $k = 1,2$:
\begin{align*}
(h_k)_\text{ghost} = (h_k)_0 + \frac{(h_k)_\text{aux} - (h_k)_0}{(x_k)_\text{aux} - (x_k)_0} \left( (x_k)_\text{ghost} - (x_k)_0 \right)
\end{align*}
where the auxiliary node on the continuous boundary is approximated by 
\begin{align*}
x_\text{aux} = \frac{1}{2} \left( (x_1)_1 + (x_1)_0 \right) &, \quad y_\text{aux} = \frac{1}{2} \left( (x_2)_2 + (x_2)_0 \right) \quad \text{and} \\
(h_k)_\text{aux} &= \frac{1}{2} \left( (h_k)_1 + (h_k)_2 \right).
\end{align*}
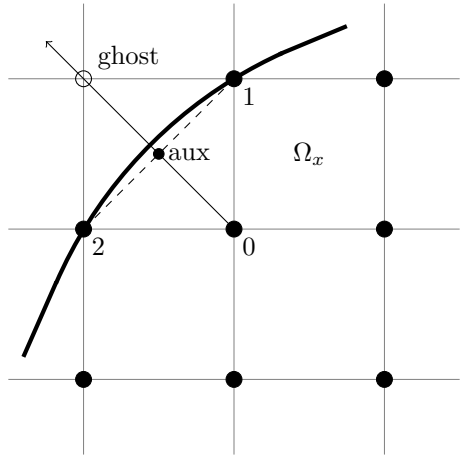
\begin{figure}[htp]
\begin{center}
	\begin{tikzpicture}
		\begin{scope}
			\draw[step=2cm,gray,very thin] (-1,-3) grid (5,3);
			\foreach \x in {0,2,4} {
				\foreach \y in {-2,0} {
					\filldraw (\x,\y) circle (3pt);
				}
			}
			\draw (0,2) circle (3pt);
			\draw (0.6,2) node [above] {ghost};
			\filldraw (2,2) circle (3pt);
			\filldraw (4,2) circle (3pt);
			\draw (2.2,2)  node [align=right, below] {1};
			\draw (0.2,0)  node [align=right, below] {2};
			\draw (2.2,0)  node [align=right, below] {0};
			\draw (3,1) node {$\Omega_x$};
			\draw [rotate around={131.5:(2,2)}] [ultra thick] (2,2) arc (-10:25:6cm);
			\draw [ultra thick] (-0.36,-0.7) -- (-0.8,-1.7);
			\draw [rotate around={131.5:(2,2)}] [ultra thick] (2,2) arc (-10:-17:6cm);
			\draw [ultra thick] (2.6,2.325) -- (3.5,2.7);
			\draw[dashed] (0,0) -- (2,2);
			\filldraw (1,1) circle (2pt) node [right] {aux};
			\draw [->] (2,0) -- (-0.5,2.5);
		\end{scope}		
	\end{tikzpicture}
 \caption{Illustration of the boundary: The nodes 1 and 2 lie directly on the continuous boundary, and their behaviour is prescribed by the Dirichlet data $\psi$. For the node 0, the stencil for the update scheme only can be applied with the help of an interpolation since the values of the ghost node are not available. The average of the values of the nodes 1 and 2 are used to create an auxiliary node which corresponds to a slightly `shifted' boundary.}
	\label{fig:interpolation}
	\end{center}
\end{figure}

We use the CFL condition 
\begin{align*}
\frac{\nu_x \Delta t}{\Delta x}  + \frac{\nu_\text{y} \Delta t}{\Delta y} \leq 1,
\end{align*} 
where $\Delta x := \min \Delta x_i$ and $\Delta y := \min \Delta y_j$, in order to determine a suitable time step $\Delta t$. The maximal propagation speeds are bounded from above by $\nu_x, \nu_\text{y} \leq \sqrt{(\lambda + 2 \mu)/\rho}$ with $\rho := \min \rho^0_{i,j} > 0$.

\section{Application in motion compensation}
\begin{figure}[htp] 
\begin{center}
\includegraphics[scale=0.7]{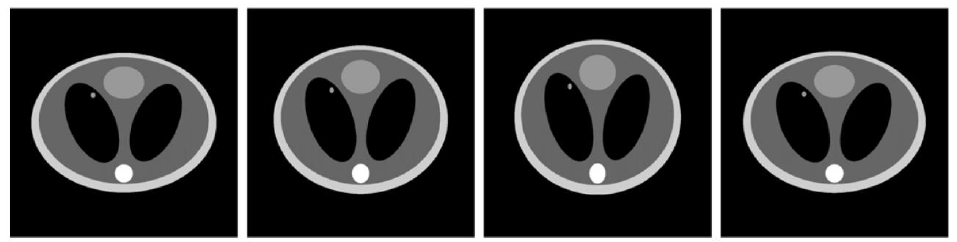}
\includegraphics[scale=0.36]{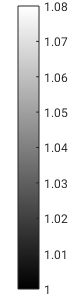}
\end{center}
\caption{Cross-section of the numerical phantom during one cycling breath. The first image corresponds to the reference state, the second and third image correspond to the body after a quarter and after one half of a breathing cycle, respectively. The fourth image illustrates the body after one period when the initial configuration is reached again.} 
\label{fig:cyclebreath}
\end{figure}

We evaluate the motion estimation approach on simulated CT data. For this purpose, we consider a thorax phantom representing a cross-section of a chest, see Figure \ref{fig:cyclebreath} left. Following from \cite{exact1}, its respiratory motion is modelled by an affine deformation, more precisely by  
\begin{align*}
\Phi(t,x) =  \begin{pmatrix}  s(t)^{-1}  & 0 \\ 0 & s(t)\end{pmatrix} \left( x - \left( \begin{array}{c}  0.44 \cdot (s(t)-1)  \\0\end{array} \right) \right)
\end{align*}
with $s(t) = 0.05 \cdot \cos(0.04 \cdot t) + 0.95$. The deformation during one breathing cycle is illustrated in the sequence of pictures in Figure \ref{fig:cyclebreath}. The phantom represents a cross-section of a simulated chest.\\

The Radon data of this dynamic object are computed for $660$ source positions, uniformly distributed over the upper half sphere, and $451$ discrete detector points uniformly distributed over $[-1,1]$ (since the support of the phantom is contained in the unit disk at all time instances). Our reconstructions and - later on - all simulations of the PDE are run on a 257x257 grid.\\

\begin{figure}[htp]
\centering
	\begin{subfigure}{0.4\textwidth}
		\includegraphics[width=\textwidth]{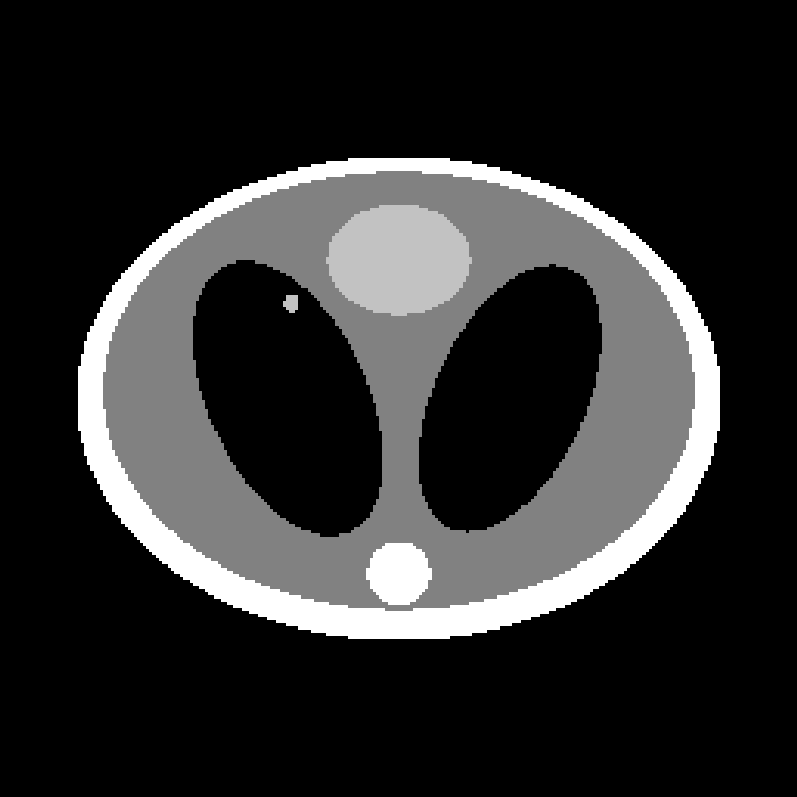}
		\subcaption{Original phantom}
	\end{subfigure}
\quad
	\begin{subfigure}{0.4\textwidth}
		\includegraphics[width=\textwidth]{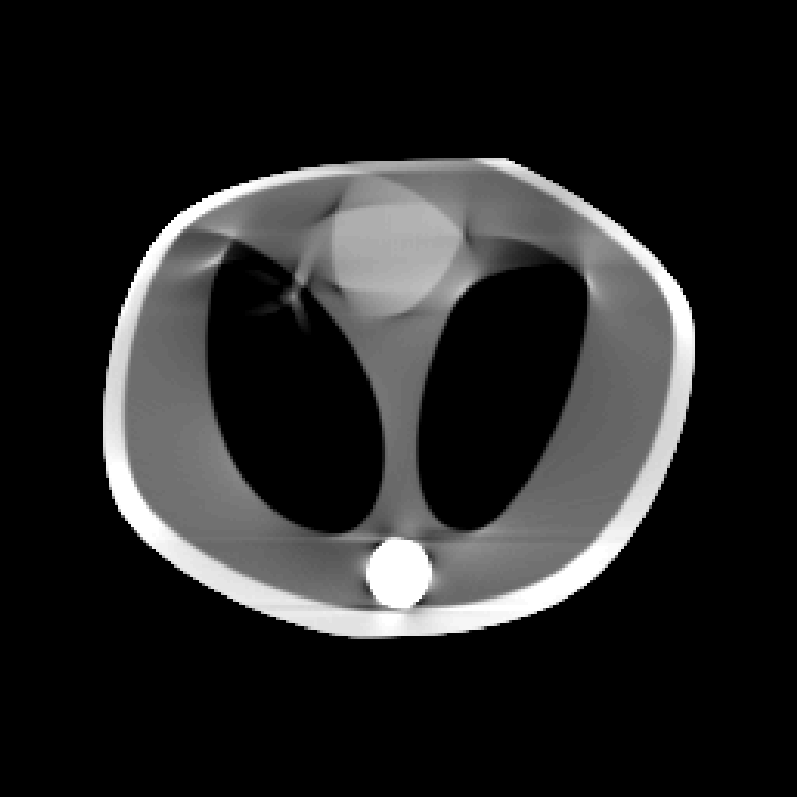} 
		\subcaption{Static reconstruction. }
	\end{subfigure}
	
		\begin{subfigure}{0.4\textwidth}
		\includegraphics[width=\textwidth]{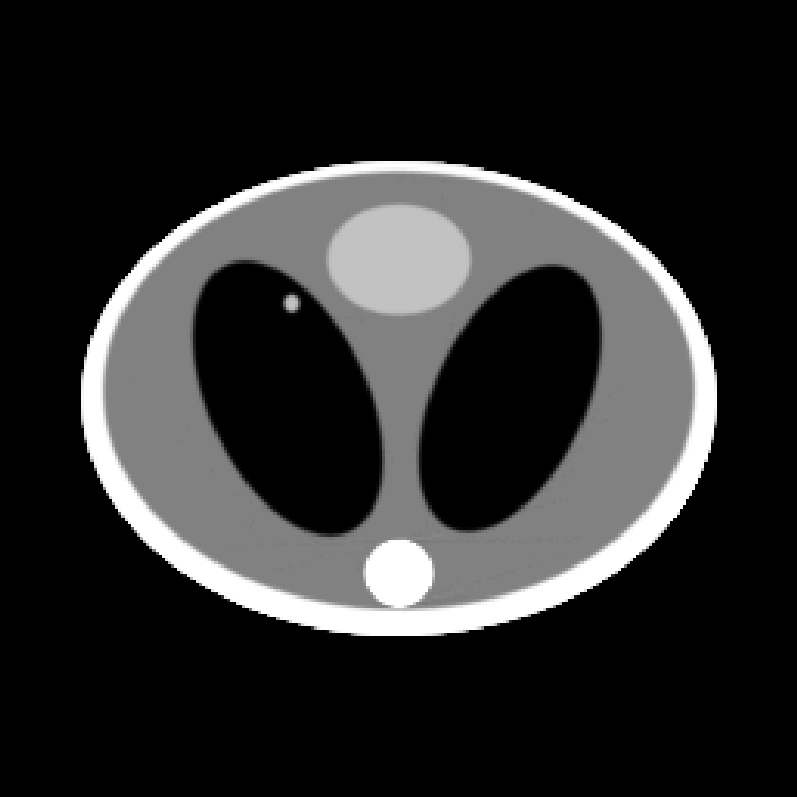}			
		\subcaption{Dynamic reconstruction with exact motion information.  \newline}
	\end{subfigure}
	\quad
	\begin{subfigure}{0.4\textwidth}
		\includegraphics[width=\textwidth]{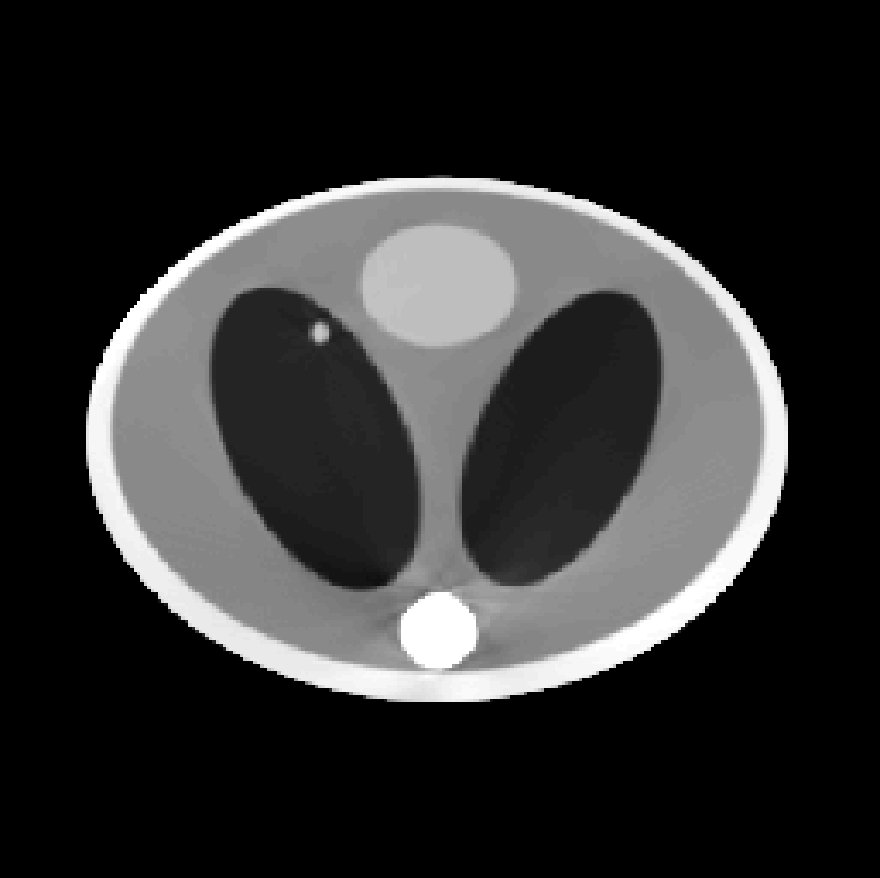}	
		\subcaption{Dynamic reconstruction with motion information from solving the PDE with analytical boundary data.}
	\end{subfigure}
	\begin{subfigure}{0.4\textwidth}
		\includegraphics[width=\textwidth]{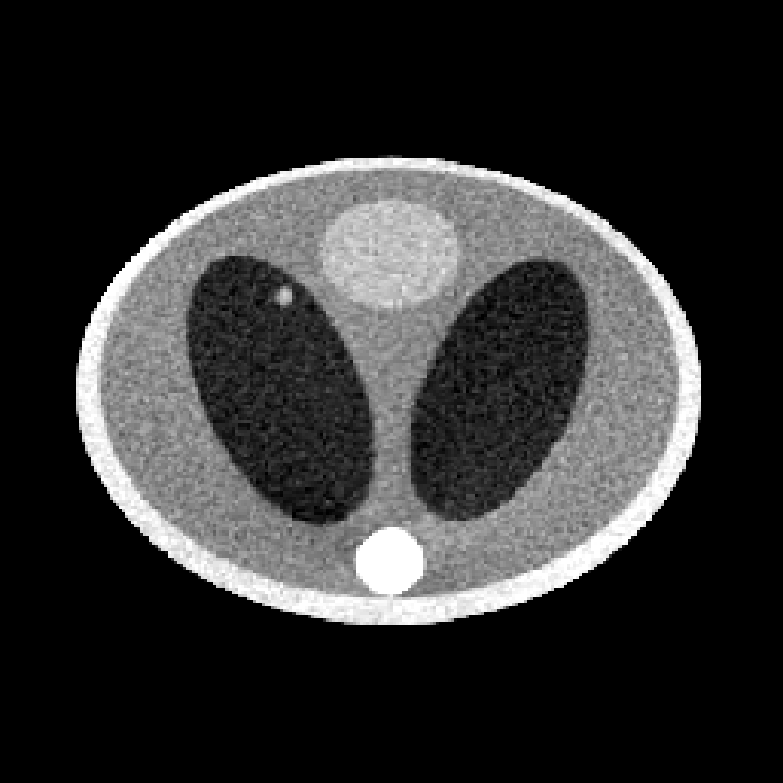}	
		\subcaption{Dynamic reconstruction from noisy Radon data with motion information from solving the PDE with analytical boundary data.}
	\end{subfigure}
	\caption{Static and dynamic reconstruction results of the initial state function.} 
	\label{fig:reconstructions-analytical}
\end{figure}

If one does not take into account that the object was moving during data acquisition and applies a static reconstruction algorithm to the dynamic data, an image of poor quality with motion artefacts such as blurring, streaking etc. is obtained, see Figure  \ref{fig:reconstructions-analytical}(B). This motivates the need for motion compensation and hence motion estimation strategies.
 
 As motion compensation algorithm, we use the strategy specified in Section \ref{Sect:CompensationCT} with the Gaussian function as low-pass filter. The result of this algorithm with exact motion information $\Phi$ is shown in Figure \ref{fig:reconstructions-analytical}(C). We observe that all components are indeed correctly reconstructed without motion artefacts, i.e. the motion is well compensated for, and in accordance to \cite{hahnkienle19}, we obtain a good approximation to the original initial state, cf. Figure \ref{fig:reconstructions-analytical}(A). However, in practice, the exact motion information is typically unknown.
 
 Thus, our goal is now to evaluate our proposed motion estimation strategy, i.e. the (discrete) deformation fields $\Phi_t$ are computed by solving the Navier-Cauchy equation with available initial and boundary data. Since the reconstruction part of the algorithm is already established in the literature, for more details about the reconstruction quality, stability regarding noise in the Radon data and the effectiveness of this part we refer to \cite{H14, hahnkienle19}. In particular, its regularizing property in order to cope with noisy measurements has been illustrated with examples from CT \cite{H14} and photoacoustic tomography \cite{hahnkienle19}. \\ 
 As discussed in Remark 5, 
 having only prescribed boundary data of a periodic and elastic movement, the Navier-Cauchy equation is a basic model to approximate the motion of the internal thorax.
 
\begin{figure}[htp]
\begin{center}
	\includegraphics[width=0.2\textwidth]{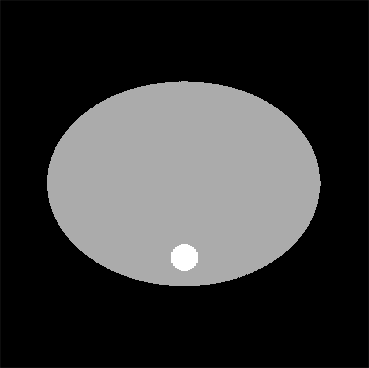}\,
	\includegraphics[width=0.2\textwidth]{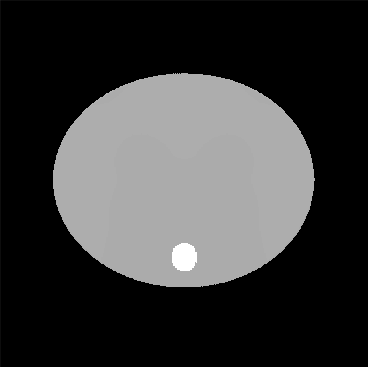}\,
	\includegraphics[width=0.2\textwidth]{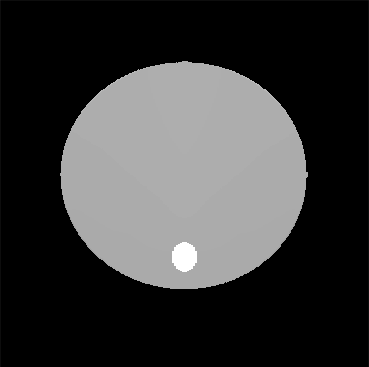}\,
	\includegraphics[width=0.2\textwidth]{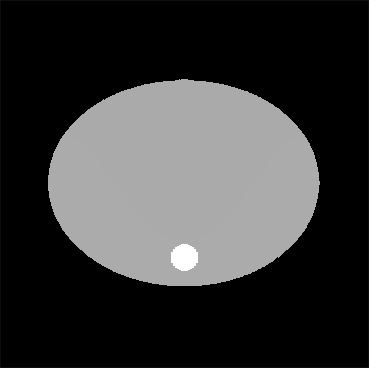}
	\caption{Illustration of the numerical solution of the Navier-Cauchy equation with analytical boundary data. The initial density distribution used for solving the Navier-Cauchy equation is given in the first image. The second, third and fourth image correspond to the configurations after a quarter, after one half and after a full period of the breathing cycle.} \label{fig:numerical-solution-PDE}
	\end{center}
\end{figure}

 First, we discuss the initial data corresponding to the initial density distribution $\hat{\rho}$. As discussed in Remark \ref{remark:initial-density}, this initial parameter is strongly linked to the searched-for initial state function $f_0$ which is why we propose to use a simplified prior instead. The one used for our simulation is shown in the first image of Figure \ref{fig:numerical-solution-PDE}.  This prior only distinguishes between spine and soft tissue, where the respective values are initialized with standard values $\hat{\rho} = 1.85\cdot 10^3$\,kg/m$^3$ for the spine and $\hat{\rho} = 1.05\cdot 10^3$\,kg/m$^3$ for the rest. This is indeed a reasonable prior in practice since the only component considered in the interior - the spine - typically does not move, so it can be extracted from a static reconstruction, cf. Figure \ref{fig:reconstructions-analytical}(B). This prior can optionally be improved by an iteration between the motion estimation with given $\hat{\rho}$ and image reconstructions, which then update $\hat{\rho}$ again.
 
Finding realistic values for the Lam\'{e}-coefficients of human tissue is a research topic by itself.  It is hard to quantify them  and they differ depending on the study \cite{Werner13}. We assume a uniform motion behavior of all (soft) tissues and restrict ourselves to one set of values for the whole thorax. The coefficients are averaged to $\lambda = 3.46$\,kPa and $\mu = 1.48$\,kPa. This simplifying assumption is reasonable for a first approach and also yields promising results in combination with the reconstruction algorithm.  \\

Regarding the boundary data, we test several configurations. First, we use the exact analytical positions of the boundary. The result for solving the respective PDE as described in Section 4 is  illustrated in Figure \ref{fig:numerical-solution-PDE}. Then, incorporating its solution as motion information in our dynamic reconstruction algorithm provides the reconstruction result shown in Figure \ref{fig:reconstructions-analytical}(D). Furthermore, in Figure \ref{fig:reconstructions-analytical}(E) a sample of noise uniformly distributed in $[-0.003, 0.003]$ was added to the Radon data in order to illustrate stability of the motion compensation algorithm. In both cases the motion of the phantom is well compensated for and the small tumour is clearly visible. This shows that determining deformation fields by solving the Navier-Cauchy equation constitutes a valuable motion estimation strategy.

\begin{figure}[htp]
	\begin{subfigure}{0.45\textwidth}
		\includegraphics[width=\textwidth]{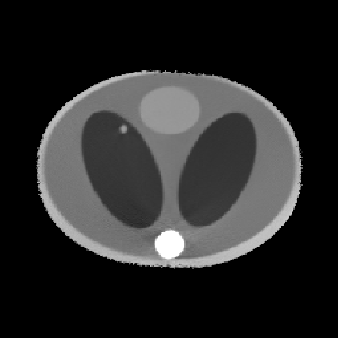}	
		\subcaption{Result for noisy boundary data with standard deviation 0.1.}
	\end{subfigure}
	\quad
	\begin{subfigure}{0.45\textwidth}
		\includegraphics[width=\textwidth]{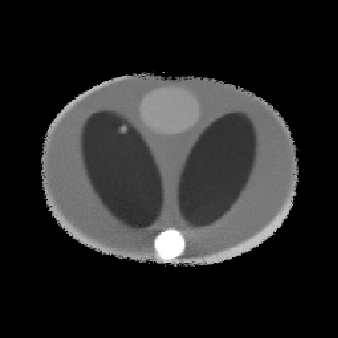}	
		\subcaption{Result for noisy boundary data with standard deviation 0.25.}
	\end{subfigure}
	\caption{Dynamic reconstruction with motion information from solving the PDE with noisy boundary data.}
	\label{fig:reconstructions-noise}
\end{figure}

In practice, the boundary positions might be determined by attaching markers at the surface of the object. If these positions are determined by measurements, they will be subject to small measurement errors. Thus, in order to test stability with respect to the boundary data, we next add a sample of noise to the (analytical) boundary positions. The noise is generated as normal distribution around 0 with standard deviation 0.1 and 0.25, respectively. In Figure \ref{fig:reconstructions-noise} we see that the reconstruction near the boundary is affected. More precisely, due to the inexact boundary positions, the boundary in the reconstruction appears fuzzy. However, the motion in the interior of the phantom is still well compensated for. All interior components, which correspond to the relevant searched-for information, including the small tumour, are still clearly recognizable, in particular in comparison to the static reconstruction, cf. Figure \ref{fig:reconstructions-analytical}(B).

\begin{figure}[htp]
	\begin{subfigure}{0.45\textwidth}
		\includegraphics[width=\textwidth]{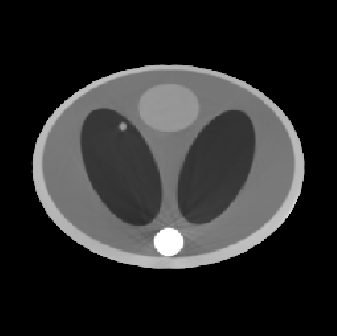}	
		\subcaption{Result for 32 prescribed boundary nodes.}
	\end{subfigure}
	\quad
	\begin{subfigure}{0.45\textwidth}
		\includegraphics[width=\textwidth]{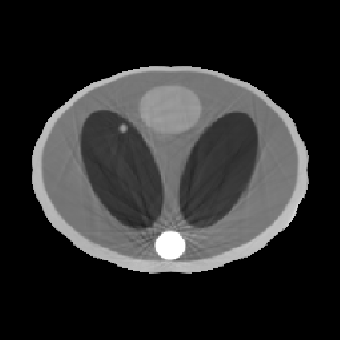}	
		\subcaption{Result for 16 prescribed boundary nodes.}
	\end{subfigure}
	\caption{Dynamic reconstruction results with motion information from solving the PDE with only a small number of boundary nodes. }
	\label{fig:reconstructions-lesspoints}
\end{figure}

Further, we test the performance of the method if only a few discrete boundary positions are given. The motivation behind this experiment is that, in practice, only a limited number of markers can be attached to the surface of the object. To this end, we prescribe only 32 (and 16, respectively) grid nodes on the boundary. Between these nodes, we apply a linear interpolation. The results are displayed in Figure \ref{fig:reconstructions-lesspoints}. We obtain some artefacts since the round shape of the thorax is replaced by a polygon due to the interpolation. However, as in the case of noisy boundary data, the deformation fields obtained by solving the PDE still provide sufficient information on the motion to compensate for it in the interior and to provide an image showing clearly all inner components including the small tumour.

\section{Conclusions and Outlook}
\label{sec:conclusions}

This article provides a proof-of-concept for a motion estimation strategy in dynamic imaging, where the Navier-Cauchy equation serves as a mathematical model for small elastic deformations.  To this end, we decoupled the tasks of motion estimation and image reconstruction, i.e. the Navier-Cauchy equation is solved prior to the reconstruction step using suitable and realistic initial and boundary data. Then the calculated deformation fields are incorporated into an analytic dynamic reconstruction algorithm. Our numerical results on a thorax phantom undergoing respiratory motion illustrate that this approach can significantly reduce motion artefacts in the respective images. In particular, we discussed available boundary data and illustrated their affect on the reconstruction result.

We illustrated the usefulness and practicability of our idea: observe the thorax's boundary, find an estimation of the thorax's motion by numerically solving a PDE, and use this approximate motion data in the reconstruction algorithm. In our numerical results, we see that the reconstruction is improved a lot compared to a static reconstruction, even by using this basic motion estimation.

In this proof-of-concept study, all simulations were run on modern desktop computers. The total computational time is the sum of time needed to solve the PDE plus to solve the inverse problem. Concerning the algorithm for the inverse problem: due to its construction, the complexity of the motion compensation algorithm is the same as for standard filtered backprojection algorithms, described for instance in \cite{Natterer86}. The implementation of the numerical scheme for the PDE was not tuned regarding efficiency, so simulation for the motion data took hours. However, in work in progress we have more elaborated techniques at hand, which significantly reduce computational costs for solving the PDE. We expect that we can then run a simulation within a few minutes on a laptop computer, this includes the time for solving the inverse problem. Thus it should be feasible to use our method in applications.

In future work, a more realistic biomechanical material law than Hooke's law will be considered. More elaborated numerical schemes will then be implemented for more specific studies, also regarding computation times. A worthwhile approach might be to minimize the distance between observed and simulated displacements in combination with solving an initial boundary value problem. Additionally, more importance will be attached to the specific behavior of different parts of the thorax. For instance, as the heart follows its own cycle, it effects the lungs' motion and its influence would also be interesting to consider.

So far, we have decoupled the suggested motion estimation and compensation approach: For estimating the deformation fields, we included a rough prior on the initial density distribution $\hat{\rho}$. This prior was then improved by incorporating the computed motion information in the image reconstruction step. The next step is to study the joint parameter identification problem, i.e. 
to address the challenging task of recovering $\hat{\rho}$ directly from \eqref{eq:NavierCauchy} with the usual boundary conditions and the data constraint $\mathcal{A}_\Phi \hat{\rho} = g$.

\section*{Acknowledgments} The first and second authors are supported by the Deutsche Forschungsgemeinschaft under grant HA 8176/1-1. The third and fourth authors want to thank Matteo  Semplice for fruitful discussions.










\end{document}